\documentclass[12pt]{article} 
\usepackage{amssymb,amsmath,latexsym}
\usepackage{amsthm}
\usepackage{graphicx}

\newtheorem{theorem}{Theorem}[section]

\newtheorem{lemma}[theorem]{Lemma}

\theoremstyle{definition}

\theoremstyle{remark}

\begin{document}
\newcommand{\beq}{\begin{equation}} \newcommand{\eeq}{\end{equation}}
\newcommand{\zz}{\mathbb{Z}}
\newcommand{\cc}{\mathbb{C}}
\newcommand{\pp}{\mathbb{P}} 
\newcommand{\nn}{\mathbb{N}}
\newcommand{\rr}{\mathbb{R}}
\newcommand{\qq}{\mathbb{Q}}
\newcommand{\bm}[1]{{\mbox{\boldmath $#1$}}}
\newcommand{\con}{\mathrm{Comp}(n)}
\newcommand{\sn}{\mathfrak{S}_n} 
\newcommand{\fsp}{\mathfrak{S}_p} 
\newcommand{\fm}{\mathfrak{m}}
\newcommand{\fs}{\mathfrak{S}}
\newcommand{\hn}{{\cal H}_n(q)}
\newcommand{\cmp}{{\cal M}(P)}
\newcommand{\st}{\,:\,} 
\newcommand{\ffq}{\mathbb{F}_q}
\newcommand{\as}{\mathrm{as}}
\newcommand{\is}{\mathrm{is}}
\newcommand{\clp}{{\cal L}(P)}
\newcommand{\maj}{\mathrm{maj}}
\newcommand{\hz}{\hat{0}}
\newcommand{\ho}{\hat{1}}
\newcommand{\covby}{\lessdot} 
\newcommand{\cov}{\gtrdot}
\newcommand{\lgn}{\mathrm{len}}
\newcommand{\dis}{\displaystyle}
\newcommand{\comaj}{\mathrm{comaj}}
\newcommand{\modd}[1]{\,(\mathrm{mod}\, #1)}
\newcommand{\be}{\begin{enumerate}}
\newcommand{\ee}{\end{enumerate}}
\newcommand{\beas}{\begin{eqnarray*}}
\newcommand{\eeas}{\end{eqnarray*}}
\newcommand{\bea}{\begin{eqnarray}}
\newcommand{\eea}{\end{eqnarray}}

\thispagestyle{empty}

\vskip 20pt
\begin{center}
{\large\bf Promotion and Evacuation}\footnote{This material is based
  upon work supported by 
  the National Science Foundation under Grant No.~0604423.}
\vskip 15pt
{\bf Richard P. Stanley}\\[.1in]
{\it Department of Mathematics, Massachusetts Institute of
Technology}\\
{\it Cambridge, MA 02139, USA}\\
{\texttt{rstan@math.mit.edu}}\\[.2in]
{\bf\small version of 19 July 2008}\\[.2in]
{Dedicated to Anders Bj\"orner on the occasion of his sixtieth
  birthday}\\ [.3in]
\end{center}

\begin{abstract}
Promotion and evacuation are bijections on the set of linear
extensions of a finite poset first defined by Sch\"utzenberger.
This paper surveys the basic properties of these two operations and
discusses some generalizations. 
\end{abstract}

\section{Introduction.} \label{sec1} 
\indent Promotion and evacuation are bijections on the set of linear
extensions of a finite poset.  Evacuation first arose in the theory of
the RSK algorithm, which associates a permutation in the symmetric
group $\sn$ with a pair of standard Young tableaux of the same shape
\cite[pp.~320--321]{ec2}.  Evacuation was described by M.-P.\
Sch\"utzenberger \cite{schutz1} in a direct way not involving the RSK
algorithm. In two follow-up papers \cite{schutz2}\cite{schutz3}
Sch\"utzenberger extended the definition of evacuation to linear
extensions of any finite poset. Evacuation is described in terms of a
simpler operation called promotion.  Sch\"utzenberger established many
fundamental properties of promotion and evacuation, including the
result that evacuation is an involution.  Sch\"utzenberger's work was
simplified by Haiman \cite{haiman} and Malvenuto and Reutenauer
\cite{ma-re}, and further work on evacuation was undertaken by a
number of researchers (discussed in more detail below).

In this paper we will survey the basic properties of promotion and
evacuation. We will then discuss some generalizations. In particular,
the linear extensions of a finite poset $P$ correspond to the maximal
chains of the distributive lattice $J(P)$ of order ideals of $P$. We
will extend promotion and evacuation to bijections on the vector space
whose basis consists of all maximal chains of a finite graded poset
$Q$. The case $Q=B_n(q)$, the lattice of subspaces of the vector space
$\ffq^n$, leads to some results on expanding a certain product in the
Hecke algebra $\hn$ of $\sn$ in terms of the standard basis $\{T_w\st
w\in\sn\}$.

I am grateful to Kyle Petersen for some helpful comments on the first
version of this paper.

\section{Basic results.} \label{sec2} 

\begin{figure}
\centering
\centerline{\includegraphics{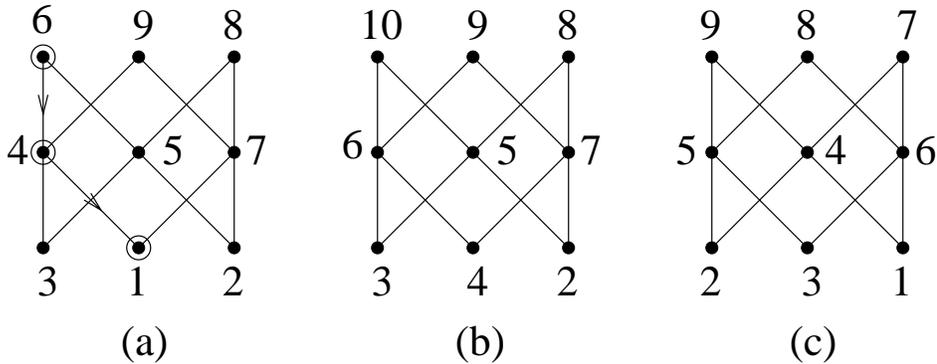}}
\caption{The promotion operator $\partial$ applied to a linear
  extension}
\label{fig:evac1}
\end{figure}

We begin with the original definitions of promotion and evacuation due
to Sch\"utzenberger. Let $P$ be a $p$-element poset. We write $s\covby
t$ if $t$ \emph{covers} $s$ in $P$, i.e., $s<t$ and no $u\in P$
satisfies $s<u<t$. The set of all linear extensions of $P$ is denoted
$\clp$. Sch\"utzenberger regards a linear extension as a bijection
$f\colon P\rightarrow[p]=\{1,2,\dots,p\}$ such that if $s<t$ in $P$,
then $f(s)<f(t)$.  (Actually, Sch\"utzenberger considers bijections
$f\colon P\rightarrow \{k+1,k+1,\dots,k+p\}$ for some $k\in\zz$, but
we slightly modify his approach by always ensuring that $k=0$.)  Think
of the element $t\in P$ as being labelled by $f(t)$. We now define a
bijection $\partial\colon \clp\rightarrow \clp$, called
\emph{promotion}, as follows. Let $t_1\in P$ satisfy $f(t_1)=1$.
Remove the label 1 from $t_1$. Among the elements of $P$ covering
$t_1$, let $t_2$ be the one with the smallest label $f(t_2)$.  Remove
this label from $t_2$ and place it at $t_1$. (Think of ``sliding'' the
label $f(t_2)$ down from $t_2$ to $t_1$.) Now among the elements of
$P$ covering $t_2$, let $t_3$ be the one with the smallest label
$f(t_3)$.  Slide this label from $t_3$ to $t_2$. Continue this process
until eventually reaching a maximal element $t_k$ of $P$.  After we
slide $f(t_k)$ to $t_{k-1}$, label $t_k$ with $p+1$. Now subtract 1
from every label. We obtain a new linear extension $f\partial \in
{\cal L}(P)$.  Note that we let $\partial$ operate on the
\emph{right}. Note also that $t_1\covby t_2 \covby \cdots \covby t_k$
is a maximal chain of $P$, called the \emph{promotion chain} of $f$.
Figure~\ref{fig:evac1}(a) shows a poset $P$ and a linear extension
$f$. The promotion chain is indicated by circled dots and arrows.
Figure~\ref{fig:evac1}(b) shows the labeling after the sliding
operations and the labeling of the last element of the promotion chain
by $p+1=10$. Figure~\ref{fig:evac1}(c) shows the linear extension
$f\partial$ obtained by subtracting 1 from the labels in
Figure~\ref{fig:evac1}(b).

It should be obvious that $\partial\colon {\cal L}(P)\rightarrow {\cal
L}(P)$ is a bijection. In fact, let $\partial^*$ denote \emph{dual
promotion}, i.e., we remove the \emph{largest} label $p$ from some
element $u_1\in P$, then slide the \emph{largest} label of an element
covered by $u_1$ up to $u_1$, etc. After reaching a minimal element
$u_k$, we label it by 0 and then add 1 to each label, obtaining
$f\partial^*$. It is easy to check that
  $$ \partial^{-1} = \partial^*. $$
\indent We next define a variant of promotion called
\emph{evacuation}. The evacuation of a linear extension $f\in
{\cal L}(P)$ is denoted $f\epsilon$ and is another linear extension of
$P$. First compute $f\partial$. Then ``freeze'' the label $p$ into
place and apply $\partial$ to what remains. In other words, let $P_1$
consist of those elements of $P$ labelled $1,2,\dots,p-1$ by $f\partial$,
and apply $\partial$ to the restriction of
$\partial f$ to $P_1$. Then freeze the label $p-1$ and apply
$\partial$ to the $p-2$ elements that remain. Continue in this way
until every element has been frozen. Let $f\epsilon$ be the linear
extension, called the \emph{evacuation} of $f$, defined by the frozen
labels.  

\medskip
\textsc{Note.} A standard Young tableau of shape $\lambda$ can be
identified in an obvious way with a linear extension of a certain
poset $P_\lambda$. Evacuation of standard Young tableaux has a nice
geometric interpretation connected with the nilpotent flag variety.
See van Leeuwen \cite[{\S}3]{leeu} and Tesler
\cite[Thm.~5.14]{tesler}. 

\begin{figure}
\centering
\centerline{\includegraphics{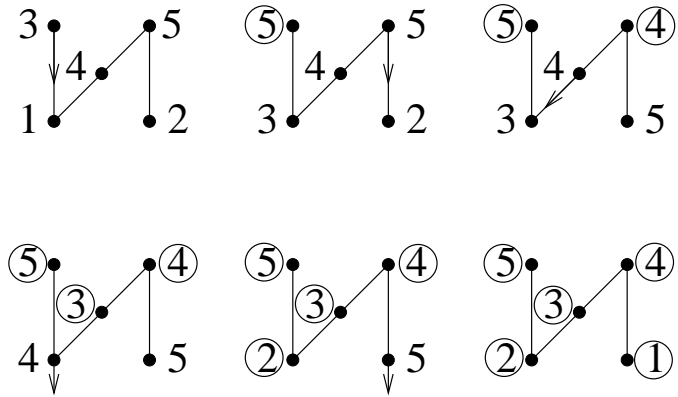}}
\caption{The evacuation of a linear extension $f$}
\label{fig:evac2}
\end{figure}

\begin{figure}
\centering
\centerline{\includegraphics{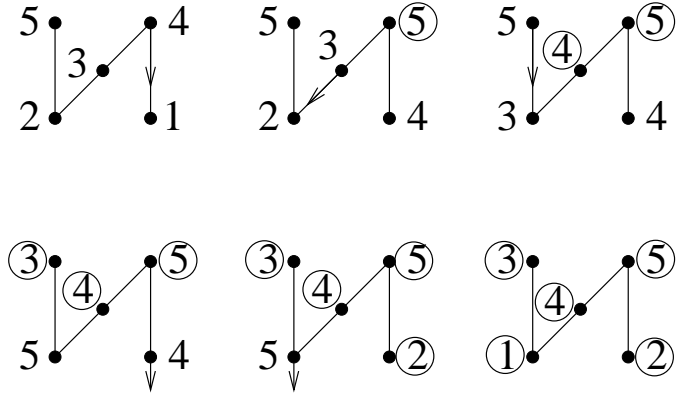}}
\caption{The linear extension evac(evac($f$))}
\label{fig:evac3}
\end{figure}

\medskip
Figure~\ref{fig:evac2} illustrates the evacuation of a linear
extension $f$. The promotion paths are shown by arrows, and the frozen
elements are circled. For ease of understanding we don't subtract 1
from the unfrozen labels since they all eventually disappear. The
labels are always frozen in descending order $p,p-1,\dots,1$. 
Figure~\ref{fig:evac3} shows the evacuation of $f\epsilon$, where $f$
is the linear extension of Figure~\ref{fig:evac2}. Note that
(seemingly) miraculously we have $f\epsilon^2=f$. This example
illustrates a fundamental property of evacuation given by
Theorem~\ref{thm:evac}(a) below.

We can define \emph{dual evacuation} analogously to dual promotion. In
symbols, if $f\in {\cal L}(P)$ then define $f^*\in{\cal L}(P^*)$ by
$f^*(t) = p+1-f(t)$. Thus
  $$ f\epsilon^* = (f^*\epsilon)^*. $$
We can now state three of the four main results obtained by
Sch\"utzenberger.  

\begin{theorem} \label{thm:evac}
Let $P$ be a $p$-element poset. Then the operators $\epsilon$,
$\epsilon^*$, and $\partial$ satisfy the following properties.
  \be\item[(a)]  Evacuation is an involution, i.e., $\epsilon^2=1$ (the
  identity operator). 
   \item[(b)]    $\partial^p=\epsilon\epsilon^*$
   \item[(c)]    $\partial\epsilon = \epsilon\partial^{-1}$
  \ee
\end{theorem}

Theorem~\ref{thm:evac} can be interpreted algebraically as follows.
The bijections $\epsilon$ and $\epsilon^*$ generate a subgroup $D_P$
of the symmetric group $\fs_{{\cal L}(P)}$ on all the linear
extensions of $P$. Since $\epsilon$ and (by duality) $\epsilon^*$ are
involutions, the group they generate is a dihedral group $D_P$
(possibly degenerate, i.e., isomorphic to $\{1\}$, $\zz/2\zz$, or
$\zz/2\zz \times \zz/2\zz$) of order 1 or $2m$ for some $m\geq 1$. If
$\epsilon$ and $\epsilon^*$ are not both trivial (which can only
happen when $P$ is a chain), so they generate a group of order $2m$,
then $m$ is the order of $\partial^p$. In general the value of $m$, or
more generally the cycle structure of $\partial^p$, is mysterious. For
a few cases in which more can be said, see Section~\ref{sec:explicit}.

The main idea of Haiman \cite[Lemma~2.7, and page~91]{haiman} (further
developed by Malvenuto and Reutenauer \cite{ma-re}) for proving
Theorem~\ref{thm:evac} is to write linear extensions as \emph{words}
rather than functions and then to describe the actions of $\partial$
and $\epsilon$ on these words. The proof then becomes a routine
algebraic computation. Let us first develop the necessary algebra in a
more general context.

Let $M$ be a monoid, i.e., a set with a binary operation (denoted by
juxtaposition) that is associative and has an identity $1$. Let
$\tau_1,\dots,\tau_{p-1}$ be elements of $M$ satisfying 
  \beq\label{eq:evacrel} \begin{array}{rl}
    \tau_i^2 =  1, & 1\leq i\leq p-1\\[.5em]  
          \tau_i \tau_j  =  \tau_j \tau_i, &\mbox{if}\ |i-j|>1. 
    \end{array} \eeq
Some readers will recognize these relations as a subset of the Coxeter
relations defining the symmetric group $\fs_p$. Define the following
elements of $M$: 
  \beas
  \delta & = & \tau_1 \tau_2\cdots \tau_{p-1}\\ \gamma\ =\
  \gamma_p & =  & \tau_1 \tau_2\cdots \tau_{p-1}\cdot 
    \tau_1 \tau_2 
    \cdots \tau_{p-2}\cdots \tau_1\tau_2\cdot \tau_1\\
   \gamma^* & = &  \tau_{p-1}\tau_{p-2}\cdots\tau_1\cdot \tau_{p-1}
    \tau_{p-2}\cdots \tau_2\cdots \tau_{p-1}\tau_{p-2}\cdot
   \tau_{p-1}. \eeas

\begin{lemma} \label{lemma:algev}
In the monoid $M$ we have the following identities:
  \be\item[\rm{(a)}]  $\gamma^2=(\gamma^*)^2=1$ 
   \item[\rm{(b)}]    $\delta^p=\gamma\gamma^*$
   \item[\rm{(c)}]    $\delta\gamma = \gamma\delta^{-1}$.
  \ee 
\end{lemma}

\proof
(a) Induction on $p$. For $p=2$, we need to show that $\tau_1^2=1$,
which is given. Now assume for $p-1$. Then
  $$ \gamma_p^2 = \tau_1\tau_2\cdots \tau_{p-1}\cdot \tau_1\cdots
  \tau_{p-2} \cdots \tau_1\tau_2\tau_3\cdot\tau_1\tau_2\cdot\tau_1 $$
 \vspace{-2em}
    $$ \qquad
   \cdot \tau_1\tau_2\cdots \tau_{p-1}\cdot \tau_1\cdots
  \tau_{p-2} \cdots \tau_1\tau_2\tau_3\cdot\tau_1\tau_2\cdot\tau_1.
  $$
We can cancel the two middle $\tau_1$'s since they appear
consecutively. We can then cancel the two middle $\tau_2$'s since they
are now consecutive. We can then move one of the middle $\tau_3$'s
past a $\tau_1$ so that the two middle $\tau_3$'s are consecutive and
can be cancelled. Now the two middle $\tau_4$'s can be moved to be
consecutive and then cancelled. Continuing in this way, we can cancel
the two middle $\tau_i$'s for all $1\leq i\leq p-1$. When this
cancellation is done, what remains is the element $\gamma_{p-1}^2$,
which is 1 by induction.

(b,c) Analogous to (a). Details are omitted.
\qed

\begin{figure}
\centering
\centerline{\includegraphics{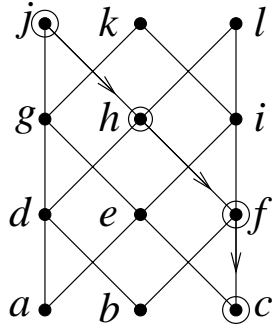}}
\caption{The promotion chain of the linear extension $cabdfeghjilk$}
\label{fig:evac4}
\end{figure}

\medskip
\emph{Proof of Theorem~\ref{thm:evac}.}  A glance at
Theorem~\ref{thm:evac} and Lemma~\ref{lemma:algev} makes it obvious
that they should be connected. To see this connection, regard the
linear extension $f\in {\cal L}(P)$ as the word (or permutation of
$P$) $f^{-1}(1),\dots,f^{-1}(p)$. For $1\leq i\leq p-1$ define
operators $\tau_i\colon {\cal L}(P)\rightarrow {\cal L}(P)$ by
  \beq \tau_i(u_1u_2\cdots u_p)=\left\{ \begin{array}{rl}
       u_1u_2\cdots u_p, & \mathrm{if}\ u_i\ \mathrm{and}\ u_{i+1}\
        \mathrm{are}\\ & \ \ \mathrm{comparable\ in}\ P\\[.2em]
       u_1u_2\cdots u_{i+1}u_i\cdots u_p, & \mathrm{otherwise}.   
           \end{array} \right. \label{eq:taui} \eeq
Clearly $\tau_i$ is a bijection, and the $\tau_i$'s satisfy the
relations \eqref{eq:evacrel}. By Lemma~\ref{lemma:algev}, the proof of
Theorem~\ref{thm:evac} follows from showing that
  $$ \partial = \delta:= \tau_1\tau_2\cdots \tau_{p-1}. $$
Note that if $f=u_1 u_2\cdots u_p$, then $f\delta$ is obtained as
follows. Let $j$ be the least integer such that $j>1$ and $u_1<u_j$.
Since $f$ is a linear extension, the elements $u_2, u_3,\dots,
u_{j-1}$ are incomparable with $u_1$. Move $u_1$ so it is between
$u_{j-1}$ and $u_j$. (Equivalently, cyclically shift the sequence $u_1
u_2\cdots u_{j-1}$ one unit to the left.) Now let $k$ be the least
integer such that $k>j$ and $u_j<u_k$. Move $u_j$ so it is between
$u_{k-1}$ and $u_k$. Continue in this way reaching the end. For
example, let $z$ be the linear extension $cabdfeghjilk$ of the poset
in Figure~\ref{fig:evac4} (which also shows the evacuation chain for
this linear extension). (We denote the linear extension for this one
example by $z$ instead of $f$ since we are denoting one of the
elements of $P$ by $f$.) We factor $z$ from left-to-right into the
longest factors for which the first element of each factor is
incomparable with the other elements of the factor:
   $$ z=(cabd)(feg)(h)(jilk). $$
Cyclically shift each factor one unit to the left to obtain
$z\delta$:
  $$ z\delta = (abdc)(egf)(h)(ilkj)=abdcegfhkilj. $$ 

\indent Now consider the process of promoting the linear extension $f$
of the previous paragraph, given as a function by $f(u_i)=i$ and as a
word by $u_1u_2\cdots u_p$. The elements $u_2,\dots,u_{j-1}$ are
incomparable with $u_1$ and thus will have their labels reduced by 1
after promotion.  The label $j$ of $u_j$ (the least element in the
linear extension $f$ greater than $u_1$) will slide down to $u_1$ and
be reduced to $j-1$.  Hence $f\partial = u_2 u_3 \cdots u_{j-1}
u_1\cdots$. Exactly analogous reasoning applies to the next step of
the promotion process, when we slide the label $k$ of $u_k$ down to
$u_j$. Hence
  $$ f\partial = u_2 u_3\cdots u_{j-1}u_1\cdot u_{j+1} u_{j+2}
      \cdots u_{k-1}u_j\cdots. $$
Continuing in this manner shows that $z\delta=z\partial$, completing
the proof of Theorem~\ref{thm:evac}.
\qed

\medskip
\textsc{Note.} The operators $\tau_i\colon \clp\rightarrow \clp$ have
the additional property that $(\tau_i\tau_{i+1})^6=1$, but we see no
way to exploit this fact. 
\medskip

\begin{figure}
\centering
\centerline{\includegraphics{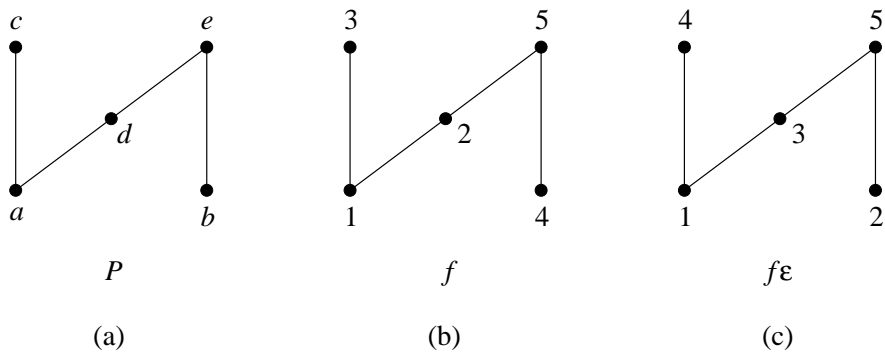}}
\caption{A poset $P$ with a linear extension and its evacuation}
\label{fig:ortr}
\end{figure}

Theorem~\ref{thm:evac} states three of the four main results of
Sch\"utzenberger. We now discuss the fourth result. Let $f\colon
P\rightarrow [p]$ be a linear extension, and apply $\partial$ $p$
times, using Sch\"utzenberger's original description of $\partial$
given at the beginning of this section. Say $f(t_1)=p$. After applying
sufficiently many $\partial$'s, the label of $t_1$ will slide down to
a new element $t_2$ and then be decreased by 1. Continuing to apply
$\partial$, the label of $t_2$ will eventually slide down to $t_3$,
etc. Eventually we will reach a minimal element $t_j$ of $P$. We call
the chain $\{t_1,t_2,\dots,t_j\}$ the \emph{principal chain} of $f$
(equivalent to Sch\"utzenberger's definition of ``orbit''), denoted
$\rho(f)$.  For instance, let $f$ be the linear extension of
Figure~\ref{fig:ortr}(b) of the poset of Figure~\ref{fig:ortr}(a).
After applying $\partial$, the label 5 of $e$ slides down to $d$ and
becomes 4. Two more applications of $\partial$ cause the label 3 to
$d$ to slide down to $a$. Thus $\rho(f)=\{a,d,e\}$.

Now apply $\partial$ to the evacuation $f\epsilon$. Let
$\sigma(f\epsilon)$ be the chain of elements of $P$ along which labels
slide, called the \emph{trajectory} of $f$. For instance,
Figure~\ref{fig:ortr}(c) shows $f\epsilon$, where $f$ is given by
Figure~\ref{fig:ortr}(b). When we apply $\partial$ to $f\epsilon$, the
label 1 of $a$ is removed, the label 3 of $d$ slides to $a$, and the
label 5 of $e$ slides to $d$. Sch\"utzenberger's fourth result is the
following.

\begin{theorem} \label{sch4}
For any finite poset $P$ and $f\in\clp$ we have
$\rho(f)=\sigma(f\epsilon)$. 
\end{theorem}

\emph{Proof} (sketch).  Regard the linear extension $f\partial^i$ of
$P$ as the word $u_{i1}u_{i2} \cdots u_{ip}$. It is not hard to check
that 
 $$ \rho(f)=\{u_{0p}, u_{1,p-1}, u_{2,p-2},\dots,u_{p-1,1}\} $$ 
(where multiple elements are counted only once). On the other hand,
let $\psi_j=\tau_1 \tau_2\cdots \tau_{p-j}$, and regard the linear
extension $f\psi_1\psi_2\cdots \psi_j$ as the word $v_{i1}v_{i2}\cdots
v_{ip}$. It is clear that $v_{ij}=u_{ij}$ if $i+j\leq p$. In
particular, $u_{i,p-i}= v_{i,p-i}$.  Moreover, $f\epsilon =
v_{2,p},v_{3,p-1},\dots,v_{p+1,1}$. We leave to the reader to check
that the elements of $\rho(f)$ written in increasing order, say
$z_1<z_2<\cdots<z_k$, form a subsequence of $f\epsilon$, since
$u_{i,p-i} = v_{i,p-i}$. Moreover, the elements of $f\epsilon$ between
$z_j$ and $z_{j+1}$ are incomparable with $z_j$. Hence when we apply
$\partial$ to $f\epsilon$, the element $z_1$ moves to the right until
reaching $z_2$, then $z_2$ moves to the right until reaching $z_3$,
etc. This is just what if means for $\sigma(f\epsilon)
=\{z_1,\dots,z_k\}$, completing the proof.  
\qed

\medskip
Promotion and evacuation can be applied to other properties of linear
extensions. We mention two such results here. For the first, let
$e(P)$ denote the number of linear extensions of the finite poset $P$.
If $A$ is the set of minimal (or maximal) elements of $P$, then it is
obvious that 
  \beq e(P) =\sum_{t\in A} e(P-t). \label{eq:eprec} \eeq
An \emph{antichain} of $P$ is a set of pairwise incomparable elements
of $P$. Edelman, Hibi, and Stanley \cite{e-h-s} use promotion to
obtain the following generalization of equation~\eqref{eq:eprec} (a
special case of an even more general theorem).

\begin{theorem} \label{thm:ehs}
Let $A$ be an antichain of $P$ that intersects every maximal
chain. Then
  $$ e(P) = \sum_{t\in A} e(P-t). $$
\end{theorem}

The second application of promotion and evacuation is to the theory of
sign balance. Fix an ordering $t_1,\dots,t_p$ of the elements of $P$,
and regard a linear extension of $f\colon P\rightarrow [p]$ as the
permutation $w$ of $P$ given by $w(t_i)=f^{-1}(i)$. A finite poset $P$
is \emph{sign balanced} if it has the same number of even linear
extensions as odd linear extensions.  It is easy to see that the
property of being sign balanced does not depend on the ordering
$t_1,\dots,t_p$. While it is difficult in general to understand the
cycle structure of the operator $\partial$ (regarded as a permutation
of the set of all linear extensions $f$ of $P$), there are situations
when we can analyze its effect on the parity of $f$. Moreover,
Theorem~\ref{thm:selfev} determines the cycle structure of $\epsilon$.
This idea leads to the following result of Stanley \cite[Cor.~2.2 and
2.4]{rs:maj}.

\begin{theorem} \label{thm:sb}
{\textrm{(a)}} Let $\#P=p$, and suppose that the length $\ell$ of
every maximal chain of $P$ satisfies $p\equiv \ell\modd{2}$. Then $P$
is sign-balanced.

{\textrm{(b)}} Suppose that for all $t\in P$, the lengths of all
  maximal chains of the principal order ideal $\Lambda_t:=\{s\in P\st
  s\leq t\}$ have the same parity. Let $\nu(t)$ denote the length of
  the longest chain of $\Lambda_t$, and set $\Gamma(P)=\sum_{t\in
    P}\nu(t)$. If $\binom p2\equiv \Gamma(P)\modd{2}$ then $P$ is
    sign-balanced. 
\end{theorem}

\section{Self-evacuation and $P$-domino tableaux} \label{sec:self}
In this section we consider self-evacuating linear extensions of a
finite poset $P$, i.e., linear extensions $f$ such that
$f\epsilon=f$. The main result asserts that the number of
self-evacuating $f\in\clp$ is equal to two other quantities associated
with $P$. We begin by defining these two other quantities.

An \emph{order ideal} of $P$ is a subset $I$ such that if $t\in I$ and
$s<t$, then $s\in I$.  A \emph{dual $P$-domino tableau} is a chain
$\emptyset=I_0\subset I_1\subset \cdots\subset I_r=P$ of order ideals
of $P$ such that $I_i-I_{i-1}$ is a two-element chain for $2\leq i
\leq r$, while $I_1$ is either a two-element or one-element chain
(depending on whether $p$ is even or odd). In particular, $r=\lceil
p/2\rceil$.

\medskip
\textsc{Note.} We use the terminology ``\emph{dual} domino tableau''
because  in \cite[{\S}4]{rs:maj} a \emph{domino tableau} is defined as
above, except that each $I_i-I_{i-1}$ is a 2-element chain for $1\leq
i\leq r-1$, and $I_r-I_{r-1}$ is either a two-element or one-element
chain. There is no difference between the two concepts if $p$ is even.

\medskip
Now assume that the vertex set of $P$ is $[p]$ and that $P$ is a
\emph{natural partial order}, i.e., if $i<j$ in $P$ then $i<j$ in
$\zz$. A linear extension of $P$ is thus a permutation $w=a_1\cdots
a_p\in \fsp$. The \emph{descent set} $D(w)$ of $w$ is defined by
  $$ D(w)=\{1\leq i\leq p-1\st a_i>a_{i+1}\}, $$
and the \emph{comajor index} $\comaj(w)$ is defined by
  \beq \comaj(w) =\sum_{i\in D(w)}(p-i). \label{eq:comaj} \eeq
(\textsc{Note.} Sometimes the comajor index is defined by $\comaj(w)=
\sum_{i\in [p-1]-D(w)}i$, but we will use equation~\eqref{eq:comaj}
here.) Set 
  $$ W'_P(x) = \sum_{w\in\clp} x^{\comaj(w)}. $$
It is known from the theory of $P$-partitions (e.g.,
\cite[{\S}4.5]{ec1}) that $W'_P(x)$ depends only on $P$ up to
isomorphism.

\medskip
\textsc{Note.} Usually in the theory of $P$-partitions one works with
the \emph{major index} maj$(w)=\sum_{i\in D(w)}i$ and with the
polynomial $W_P(x)=\sum_{w\in\clp}x^{\mathrm{maj}(w)}$. Note that if
$p$ is even then $\comaj(w)\equiv \mathrm{maj}(w)\modd{2}$, so
$W_P(-1) = W'_P(-1)$.

\begin{theorem} \label{thm:selfev}
Let $P$ be a finite natural partial order. Then the following three
quantities are equal.
  \be \item[\textrm{(i)}] $W'_P(-1)$.
  \item[\textrm{(ii)}] The number of dual $P$-domino tableaux.
  \item[\textrm{(iii)}] The number of self-evacuating linear
    extensions of $P$.
  \ee
\end{theorem}

In order to prove Theorem~\ref{thm:selfev}, we need one further result
about the elements $\tau_i$ of equation~\eqref{eq:evacrel}. 

\begin{lemma} \label{lemma:monoidse}
Let $M$ be the monoid of Lemma~\ref{lemma:algev}. Write 
  \beas \delta_i & = & \tau_1\tau_2\cdots \tau_i\\
        \delta^*_i & = & \tau_i\tau_{i-1}\cdots \tau_1. \eeas
Let $u,v\in M$. The following two conditions are equivalent.
  \be\item[(i)] 
   $$ u\delta^*_1\delta^*_3\cdots\delta^*_{2j-1} =
     v\delta^*_1\delta^*_3\cdots\delta^*_{2j-1}\cdot\delta_{2j-1}
     \delta_{2j-2}\cdots \delta_2\delta_1. $$
   \item[(ii)] $u\tau_1\tau_3\cdots\tau_{2j-1}=v$.
   \ee 
\end{lemma}

\emph{Proof of Lemma~\ref{lemma:monoidse}.}  The proof is a
straightforward extension of an argument essentially due to van
Leeuwen \cite[{\S}2.3]{vlee} and more explicitly to Berenstein and
Kirillov \cite{be-ki}. (About the same time as van Leeuwen, a special
case was proved by Stembridge \cite{stemb} using representation
theory. Both Stembridge and Berenstein-Kirillov deal with
\emph{semistandard tableaux}, while here we consider only the special
case of standard tableaux. While standard tableaux have a natural
generalization to linear extensions of any finite poset, it is unclear
how to generalize semistandard tableaux analogously so that the
results of Stembridge and Berenstein-Kirillov continue to hold.)
Induction on $j$. The case $j=1$ asserts that $u\tau_1=v\tau_1\tau_1$
if and only if $u\tau_1=v$, which is immediate from $\tau_1^2=1$. Now
assume for $j-1$, and suppose that (i) holds. First cancel
$\delta^*_{2j-1}\delta_{2j-1}$ from the right-hand side. Now take the
first factor $\tau_i$ from each factor $\delta_1,\dots,\delta_{2j-2}$
on the right-hand side and move it as far to the right as possible.
The right-hand side will then end in $\tau_{2j-2}\tau_{2j-3}\cdots
\tau_1=\delta^*_{2j-2}$. The left-hand side ends in
$\delta^*_{2j-1}=\tau_{2j-1}\delta^*_{2j-2}$.  Hence we can cancel the
suffix $\delta^*_{2j-2}$ from both sides, obtaining \beq
u\delta^*_1\delta^*_3\cdots\delta^*_{2j-3}\tau_{2j-1} =
v\delta^*_1\delta^*_3\cdots\delta^*_{2j-3}\cdot\delta_{2j-3}
\delta_{2j-4}\cdots \delta_2\delta_1. \label{eq:indeq} \eeq We can now
move the rightmost factor $\tau_{2j-1}$ on the left-hand side of
equation~\eqref{eq:indeq} directly to the right of $u$. Applying the
induction hypothesis with $u$ replaced by $u\tau_{2j-1}$ yields (ii).
The steps are reversible, so (ii) implies (i).  \qed

\medskip
\emph{Proof of Theorem~\ref{thm:selfev}.}
The equivalence of (i) and (ii) appears (in dual form) in
\cite[Theorem~5.1(a)]{rs:maj}. Namely, let $w= a_1\cdots a_p\in{\cal
  L}(P)$. Let $i$ be the least nonnegative integer (if it exists) for
which
   $$ w':=a_1\cdots a_{p-2i-2}a_{p-2i}a_{p-2i-1}a_{p-2i+1}\cdots
     a_p\in{\cal L}(P). $$ 
Note that $w''=w$. Now exactly one of $w$ and $w'$ has the descent
$p-2i-1$. The only other differences in the descent sets of $w$ and
$w'$ occur (possibly) for the numbers $p-2i-2$ and $p-2i$. Hence
$(-1)^{\comaj(w)}+(-1)^{\comaj(w')}=0$. The surviving permutations
$w=b_1\cdots b_p$ in ${\cal L}(P)$ are exactly those for which the
chain of order ideals 
  $$ \emptyset\subset \cdots \subset \{b_1,b_2,\dots,b_{p-4}\}
     \subset \{b_1,b_2,\dots,b_{p-2}\}\subset
      \{b_1,b_2,\dots,b_p\}=P $$
  is a dual $P$-domino tableau. We call $w$ a \emph{dual domino linear
    extension}; they are in bijection with dual domino tableaux. Such
  permutations $w$ can only have descents in positions $p-j$ where $j$
  is even, so $(-1)^{\comaj(w)}=1$.  Hence (i) and (ii) are equal.

To prove that (ii) and (iii) are equal, let $\tau_i$ be the operator
on ${\cal L}(P)$ defined by equation~\eqref{eq:taui}. Thus $w$ is
self-evacuating if and only if
    $$ w = w\tau_1\tau_2\cdots \tau_{p-1}\cdot \tau_1\cdots
      \tau_{p-2} \cdots
      \tau_1\tau_2\tau_3\cdot\tau_1\tau_2\cdot\tau_1. $$
On the other hand, note that $w$ is a dual domino linear extension if
and only if
   $$ w\tau_{p-1}\tau_{p-3}\tau_{p-5}\cdots \tau_h=w, $$
where $h=1$ if $p$ is even, and $h=2$ if $p$ is odd. It follows from
Lemma~\ref{lemma:monoidse} (letting $u=v=w$) that
$w$ is a dual domino linear extension if and only if 
    $$ \widetilde{w}:=w\tau_1\cdot \tau_3\tau_2\tau_1\cdot
      \tau_5\tau_4\tau_3\tau_2\tau_1\cdots \tau_m\tau_{m-1}\cdots
      \tau_1 $$
is self-evacuating, where $m=p-1$ if $p$ is even, and $m=p-2$ if $p$
is odd. The proof follows since the map
$w\mapsto\widetilde{w}$ is then a bijection between dual domino linear
extensions and self-evacuation linear extensions of $P$. 
\qed

%

\medskip
The equivalence of (i) and (iii) above is an instance of Stembridge's
``$q=-1$ phenomenon.'' Namely, suppose that an involution $\iota$ acts
on a finite set $S$. Let $f:S\rightarrow \zz$. (Usually $f$ will be a
``natural'' combinatorial or algebraic statistic on $S$.)  Then we say
that the triple $(S,\iota,f)$ exhibits the $q=-1$ phenomenon if the
number of fixed points of $\iota$ is given by $\sum_{t\in
  S}(-1)^{f(t)}$. See Stembridge
\cite{stemb2}\cite{stemb3}\cite{stemb}. The $q=-1$ phenomenon has been
generalized to the action of cyclic groups by V. Reiner, D.  Stanton,
and D. White \cite{r-s-w}, where it is called the ``cyclic sieving
phenomenon.'' For further examples of the cyclic sieving phenomenon,
see C. Bessis and V.  Reiner \cite{be-re}, H. Barcelo, D. Stanton, and
V. Reiner \cite{b-s-r}, and B. Rhoades \cite{rhoades}. In the next
section we state a deep example of the cyclic sieving phenomenon, due
to Rhoades, applied to the operator $\partial$ when $P$ is the product
of two chains.

\section{Special cases.} \label{sec:explicit} 
There are a few ``nontrivial'' classes of posets $P$ known for which
the operation $\partial^p=\epsilon\epsilon^*$ can be described in a
simple explicit way, so in particular the order of the dihedral group
$D_P$ generated by $\epsilon$ and $\epsilon^*$ can be determined.
There are also some ``trivial'' classes, such as hook shapes (a
disjoint union of two chains with a $\hat{0}$ adjoined), where it is
straightforward to compute the order of $\partial$ and $D_P$. The
nontrivial classes of posets are all connected with the theory of
standard Young tableaux or shifted tableaux, whose definition we
assume is known to the reader. A standard Young tableau of shape
$\lambda$ corresponds to a linear extension of a certain poset
$P_\lambda$ in an obvious way, and similarly for a standard shifted
tableau. (As mentioned in the introduction, Sch\"utzenberger
originally defined evacuation for standard Young tableaux before
extending it to linear extensions of any finite poset.) We will simply
state the known results here. The posets will be defined by examples
which should make the general definition clear.

\begin{figure}
\centering
\centerline{\includegraphics{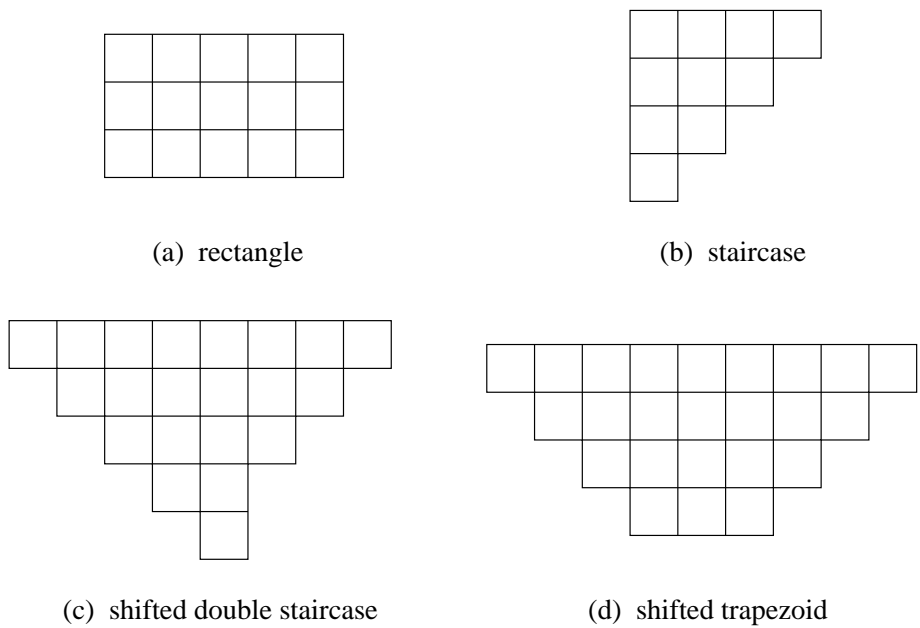}}
\caption{Some shapes and shifted shapes}
\label{fig:special}
\end{figure}

\begin{theorem} \label{thm:special}
For the following shapes and shifted shapes $P$ with a total of
$p=\#P$ squares, we have the indicated properties of $\partial^p$
and $D_P$.
  \be\item[\textrm{(a)}] Rectangles
  (Figure~\ref{fig:special}(a)). Then $f\partial^p=f$ and $D_P\cong
  \zz/2\zz$ (if $m,n>1$). Moreover, if $f=(a_{ij})$ (where we are
  regarding a linear extension of the rectangle $P$ as a labeling of
  the squares of $P$), then $f\epsilon = (p+1-a_{m+1-i,n+1-j})$. 
  \item[\textrm{(b)}] Staircases (Figure~\ref{fig:special}(b)). Then
    $f\partial^p= f^t$ (the transpose of $f$) and $D\cong \zz/2\zz
    \times \zz/2\zz$. 
  \item[\textrm{(c)}] Shifted double staircases
    (Figure~\ref{fig:special}(c)). Then $f\partial^p=f$ and $D_P\cong
  \zz/2\zz$. 
  \item[\textrm{(d)}] Shifted trapezoids
    (Figure~\ref{fig:special}(d)). Then $f\partial^p=f$ and $D_P\cong
  \zz/2\zz$. 
 \ee
\end{theorem}

Theorem~\ref{thm:special}(a) follows easily from basic properties of
\emph{jeu de taquin} due to Sch\"utzenberger \cite{schutz4} (see also
\cite[Ch.~7, Appendix~1]{ec2}) and is often attributed to
Sch\"utzenberger. We are unaware, however, of an explicit statement in
the work of Sch\"utzenberger. Part~(b) is due to Edelman and Greene
\cite[Cor.~7.23]{e-g}. Parts~(c) and (d) are due to Haiman
\cite[Thm.~4.4]{haiman}, who gives a unified approach also including
(a) and (b).

The equivalence of (i) and (iii) in Theorem~\ref{thm:selfev} was given
a deep generalization by Rhoades \cite{rhoades} when $P$ is an
$m\times n$ rectangular shape (so $p=mn$), as mentioned in the
previous section. By Theorem~\ref{thm:special}(a) we have
$f\partial^p=f$ when $P$ is a rectangular shape of size $p$. Thus
every cycle of $\partial$, regarded as a permutation of the set
$\clp$, has length $d$ dividing $p$. We can ask more generally for the
precise cycle structure of $\partial$, i.e., the number of cycles of
each length $d|n$. Equivalently, for any $d\in\zz$ (or just any
$d|p$) we can ask for the quantity
  $$ e_d(P) = \#\{f\in\clp\st f=f\partial^d\}. $$
To answer this question, define the \emph{major index} of the linear
extension $f\in\clp$ by
  $$ \mathrm{maj}(f)=\sum_i i, $$
where $i$ ranges over all entries of $P$ for which $i+1$ appears in a
lower row than $i$ \cite[p.~363]{ec2}. For instance, if $f$ is given
by
  $$ f = \begin{array}{cccc} 1 & 3 & 4 & 8\\ 2 & 5 & 6 & 11\\
            7 & 9 & 10 & 12 \end{array}, $$
then maj$(f)=1+4+6+11=22$. Let 
  $$ F(q) = \sum_{f\in\clp} q^{\maj(f)}. $$
It is well known \cite[Cor.~7.21.5]{ec2} that
  $$ F(q) = \frac{q^{n\binom m2}(1-q)(1-q^2)\cdots(1-q^p)}
                   {\prod_{t\in P}(1-q^{h(t)})}, $$
where $h(t)$ is the hook length of $t$. If say $m\leq n$, then we have
more explicitly
  $$ \hspace{-2in} \prod_{t\in P}(1-q^{h(t)}) $$
  $$ = [1][2]^2[3]^3\cdots[m]^m[m+1]^m
  \cdots [n]^m[n+1]^{m-1}[n+2]^{m-2}\cdots[n+m-1], $$
where $[i]=1-q^i$. The beautiful result of Rhoades is the following.

\begin{theorem} \label{thm:rhoades}
Let $P$ be an $m\times n$ rectangular shape. Set $p=mn$ and $\zeta =
e^{2\pi i/p}$. Then for any $d\in \zz$ we have
  $$ e_d(P) = F(\zeta^d). $$
\end{theorem}
  
Rhoades' proof of this theorem uses Kazhdan-Lusztig theory and a
characterization of the dual canonical basis of
$\cc[x_{11},\dots,x_{nn}]$ due to Skandera \cite{skan}. Several
questions are suggested by Theorems~\ref{thm:special} and
\ref{thm:rhoades}.

\be
 \item Is there a more elementary proof of Theorem~\ref{thm:rhoades}?
   For the special case of $2\times n$ and $3\times n$ rectangles, see
   \cite{p-p-r}. The authors of this paper are currently extending
   this proof to general rectangles. 
 \item Can Theorem~\ref{thm:rhoades} be extended to more general
   posets, in particular, the posets of
   Theorem~\ref{thm:special}(b,c,d)? 
 \item Can Theorem~\ref{thm:special} itself be extended to other
   classes of posets? A possible place to look is among the
   $d$-complete posets of Proctor \cite{proctor}. Some work along
   these lines is being done by Kevin Dilks (in progress at the time
   of this writing). 
\ee

\section{Generalizations.} The basic properties of evacuation given in
Sections~\ref{sec2} and \ref{sec:self} depend only on the formal
properties of the monoid $M$ defined by
equation~\eqref{eq:evacrel}. It is easy to find other examples of
operators satisfying these conditions that are more general than the
operators $\tau_i$ operating on linear extensions of posets. Hence the
theory of promotion and evacuation extends to these more general
situations. 

Let $J(P)$ denote the set of all order ideals of the finite poset $P$,
ordered by inclusion. By a well-known theorem of Birkhoff (see
\cite[Thm.~3.41]{ec1}), the posets $J(P)$ coincide with the finite
distributive lattices. There is a simple bijection \cite[{\S}3.5]{ec1}
between maximal chains $\emptyset=I_0\subset I_1\subset\cdots\subset
I_p=P$ of $J(P)$ and linear extensions of $P$, viz., associate with
this chain the linear extension $f\colon P\rightarrow [p]$ defined by
$f(t)=i$ if $t\in I_i-I_{i-1}$. In terms of the maximal chain
$\mathfrak{m}\colon \emptyset=I_0\subset I_1\subset\cdots\subset
I_p=P$ of $J(P)$, the operator $\tau_i$ on linear extensions of $P$
can be defined as follows. The interval $[I_{i-1},I_{i+1}]$ contains
either three or four elements, i.e., either $I_i$ is the unique
element satisfying $I_{i-1}\subset I_i\subset I_{i+1}$ or there is
exactly one other such element $I'$. In the former case define
$\tau_i(\fm)=\fm$; in the latter case, $\tau_i(\fm)$ is obtained from
$\fm$ by replacing $I_i$ with $I'$.

The exact same definition of $\tau_i$ can be made for any finite
graded poset, say for convenience with a unique minimal element $\hz$
and unique maximal element $\ho$, for which every interval of rank 2
contains either three or four elements. Let us call such posets
\emph{slender}. Clearly the $\tau_i$'s satisfy the conditions
\eqref{eq:evacrel}. Thus Lemma~\ref{lemma:algev} applies to the
operators $\gamma$, $\gamma^*$, and $\delta$.  (These observations
seem first to have been made by van Leeuwen \cite[{\S}2]{vlee}, after
similar results by Malvenuto and Reutenauer \cite{ma-re} in the context
of graphs rather than posets.) We also have an analogue for slender
posets $Q$ of the equivalence of (ii) and (iii) in
Lemma~\ref{lemma:monoidse}. The role of dual $P$-domino tableau is
played by \emph{dual domino chains} of $Q$, i.e., chains
$\hz=t_0<t_1<\cdots<t_r=\ho$ in $P$ for which the interval
$[t_{i-1},t_i]$ is a two-element chain for $2\leq i \leq r$, while
$[t_0,t_1]$ is either a two-element or one-element chain (depending on
whether the rank of $Q$ is even or odd). We then have that the number
of self-evacuating maximal chains of $Q$ is equal to the number of
dual domino chains of $Q$.

Some example of slender posets are Eulerian posets
\cite[{\S}3.14]{ec1}, which include face posets of regular CW-spheres
\cite{bjorner} and intervals in the Bruhat order of Coxeter groups $W$
(including the full Bruhat order of $W$ when $W$ is finite). Eulerian
posets $Q$ have the property that every interval of rank 2 contains
four elements.  Hence there are no dual domino chains when
rank$(Q)>1$, and therefore also no self-evacuating maximal chains.
Non-Eulerian slender posets include the weak order of a finite Coxeter
group \cite{bj:weak}\cite[Ch.~3]{bj-br} and face posets of regular
CW-balls. We have not systematically investigated whether there are
examples for which more can be said, e.g., an explicit description of
evacuation or the determination of the order of the dihedral group
generated by $\gamma$ and $\gamma^*$.

There is a simple example that can be made more explicit, namely, the
face lattice $L_n$ of an $n$-dimensional cross-polytope ${\cal C}_n$
(the dual to an $n$-cube). The vertices of ${\cal C}_n$ can be
labelled $1,\bar{1},2,\bar{2},\dots,n,\bar{n}$ so that vertices $i$
and $\bar{i}$ are antipodal for all $i$. A maximal chain
$\hz=t_0<t_1<\cdots<t_{n+1}=\ho$ of $L_n$ can then be encoded as a
\emph{signed permutation} $a_1,\dots,a_n$, i.e., take a permutation
$b_1,\dots,b_n$ and place bars above some subset of the $b_i$'s. 
Thus $a_i$ is the unique vertex of the face $t_i$ that does not lie in
$t_{i-1}$. Write $'$ for the reversal of the bar, i.e., $i'=\bar{i}$
and $\bar{i}'=i$. Let $w=a_1, \dots,a_n$ be a signed permutation of
$1,2,\dots, n$. Then it is easy to compute that
  \beas w\delta & = & a_2, a_3, \dots, a_n,a'_1\\
        w\gamma & = & a'_1, a_n, a_{n-1},\dots,a_2\\
        w\gamma^* & = & a'_n, a'_{n-1},\dots,a'_1\\
       w\delta^{n+1}\ =\ w\gamma\gamma^* & = & 
          a'_2, a'_3,\dots,a'_n,a_1. \eeas
Thus $\gamma\gamma^*$ has order $n$ if $n$ is odd and $2n$ if $n$ is
even. The dihedral group generated by $\gamma$ and $\gamma^*$ has
order $2n$ if $n$ is odd and $4n$ if $n$ is even.

Can the concepts of promotion and evacuation be extended to posets
that are not slender? We discuss one way to do this. Let $P$ be a
graded poset of rank $n$ with $\hz$ and $\ho$. If $\fm\colon \hz=t_0 <
t_1 < \cdots <t_n=\ho$ is a maximal chain of $P$, then we would like
to define $\fm\tau_i$ so that (1) $\tau_i^2=1$, and (2) the action of
$\tau_i$ is ``local'' at rank $i$, i.e., $\fm\tau_i$ should only
involve maximal chains that agree with $\fm$ except possibly at
$t_i$. There is no ``natural'' choice of a single chain
$\fm'=\fm\tau_i$, so we should be unbiased
and choose a linear combination of chains. Thus let $K$ be a field of
characteristic 0. Write $\cmp$ for the set of
maximal chains of $P$ and $K\cmp$ for the $K$-vector space with basis
$\cmp$. For $1\leq i\leq n-1$ define a linear operator 
$\tau_i\colon K\cmp\rightarrow K\cmp$ as follows. Let $N_i(\fm)$ be the
set of maximal chains $\fm'$ of $P$ that differ from $\fm$ exactly at
$t_i$, i.e., $\fm'$ has the form 
   $$ \fm'\colon\hz=t_0<t_1<\cdots<t_{i-1}<t'_i<t_{i+1}<\cdots<
       t_n=\ho, $$
where $t'_i\neq t_i$. Suppose that $\#N_i(\fm)=q\geq 1$. Then set
  \beq \fm\tau_i = \frac{1}{q+1}\left((q-1)\fm-2\sum_{\fm'\in
    N_i(\fm)}\fm'\right). \label{eq:tauidef} \eeq
When $q=0$ we set $\fm\tau_i=\fm$, though it would make no difference
to set $\fm\tau_i=-\fm$ to remain consistent with
equation~\eqref{eq:tauidef}.  It is easy to check that $\tau_i^2=1$.
In fact, $\pm\tau_i$ are the unique involutions of the form
$a\fm+b\sum_{\fm'\in N_i(\fm)}\fm'$ for some $a,b\in K$ with $b\neq
0$ when $q\geq 1$. It is clear that also $\tau_i \tau_j=\tau_j\tau_i$
if $|j-i|\geq 2$, so the $\tau_i$'s satisfy \eqref{eq:evacrel}. Hence
we can define promotion and evacuation on the maximal chains of any
finite graded poset so that Lemma~\ref{lemma:algev} holds, as well as
an evident analogue of the equivalence of (ii) and (iii) in
Theorem~\ref{thm:selfev}. 

The obvious question then arises: are there interesting examples? We
will discuss one example here, namely, the lattice $B_n(q)$ of
subspaces of the $n$-dimensional vector space $\ffq^n$ (ordered by
inclusion). This lattice is the ``$q$-analogue'' of the boolean
algebra $B_n$ of all subsets of the set $\{1,2,\dots,n\}$, ordered by
inclusion. The boolean algebra $B_n$ is the lattice of order ideals of
an $n$-element antichain. Hence promotion and evacuation on the
maximal chains of $B_n$ are equivalent to ``classical'' promotion and
evacuation on an $n$-element antichain $[n]=\{1,2,\dots,n\}$.  The
linear extensions of $A$ are just \emph{all} the permutations $w$ of
$[n]$, and the evacuation $w\epsilon$ of $w=a_1 a_2\cdots a_n$ is just
the reversal $a_n\cdots a_2 a_1$. Thus we are asking for a kind of
$q$-analogue of reversing a permutation.

This problem can be reduced to a computation in the Hecke algebra
$\hn$ of the symmetric group $\sn$ over the field $K$ (of
characteristic 0). Recall
(e.g., \cite[{\S}7.4]{hum}) that $\hn$ has generators $T_1,\dots,T_{n-1}$
and relations
    \beas (T_i+1)(T_i-q) & = & 0\\
      T_i T_j & = & T_j T_i,\ \ |i-j|\geq 2\\
      T_i T_{i+1} T_i & = & T_{i+1} T_i T_{i+1}. \eeas
If $q=1$ then we have $T_i^2=1$, and the above relations are just the
Coxeter relations for the group algebra $K\sn$. For $1\leq i\leq
n-1$ let $s_i$ denote the adjacent transposition $(i,i+1)\in\sn$. A
\emph{reduced decomposition} of an element $w\in\sn$ is a sequence
$(a_1,\dots,a_r)$ of integers $1\leq a_i\leq n-1$ such that
$w=s_{a_1}\cdots s_{a_r}$ and $r$ is as small as possible, namely, $r$
is the number of inversions of $w$. Define $T_w=T_{a_1}\cdots
T_{a_r}$. A standard fact about $\hn$ is that $T_w$ is independent of
the choice of reduced decomposition of $w$, and the $T_w$'s for
$w\in\sn$ form a $K$-basis for $\hn$. 


Let End$(K\cmp)$ be the set of all linear transformations $K\cmp
\rightarrow K\cmp$.  It is easy to check that the map $T_i\mapsto
\tau_i$ extends to an algebra homomorphism (i.e., a representation of
$\hn$) $\varphi\colon\hn\rightarrow \mathrm{End}(K\cmp)$. Moreover,
$\varphi$ is injective. If we fix a maximal chain $\fm_0$, then the set
$\cmp$ has a \emph{Bruhat decomposition} \cite[{\S}23.4]{f-h}
  $$ \cmp = \bigsqcup_{w\in\sn} \Omega_w, $$
where $\bigsqcup$ denotes disjoint union and $\Omega_{\mathrm{id}}=
\{\fm_0\}$. We then have
   $$ \tau_w(\fm_0) = \sum_{\fm\in\Omega_w} \fm. $$
Let $E_i=\frac{1}{q+1}(q-1-2T_i)\in\hn$, so $E_i^2=1$. It follows that 
  $$ \fm_0\epsilon = \sum_{w\in\sn} c_w(q)\sum_{\fm\in\Omega_w}w, $$
where $c_w(q)$ is defined by the Hecke algebra expansion
  \beq E_1E_2\cdots E_{n-1}E_1E_2\cdots E_{n-2}\cdots E_1E_2E_1 =
      \sum_{w\in\sn} c_w(q)T_w. \label{eq:cwdef} \eeq
For instance, when $w\in\fs_4$ the values of $c_w(q)$ are
given by
  \beas c_{1234}(q) & = & (q-1)^2/(q+1)^2\\
        c_{1243}(q) & = & -2(q-1)^3/(q+1)^4\\       
        c_{1324}(q) & = & -16q(q-1)(q^2+1)/(q+1)^6\\       
        c_{1342}(q) & = & 4(q-1)^2/(q+1)^4\\       
        c_{1423}(q) & = & 4(q-1)^2/(q+1)^4\\       
        c_{1432}(q) & = & -8(q-1)^3/(q+1)^6\\       
        c_{2134}(q) & = & -2(q-1)^3/(q+1)^4\\       
        c_{2143}(q) & = & 4(q-1)^2/(q+1)^4\\       
        c_{2314}(q) & = & -4(q-1)^4/(q+1)^6\\       
        c_{2341}(q) & = & -8(q-1)/(q+1)^4\\ 
        c_{2413}(q) & = & 0\\ 
        c_{2431}(q) & = & 16(q-1)^2/(q+1)^6\\       
        c_{3124}(q) & = & -4(q-1)^4/(q+1)^6\\       
        c_{3142}(q) & = & 0\\       
        c_{3214}(q) & = & 8(q-1)^3/(q+1)^6\\       
        c_{3241}(q) & = & 0\\       
        c_{3412}(q) & = & 16(q-1)^2/(q+1)^6\\       
        c_{3421}(q) & = & -32(q-1)/(q+1)^6\\       
        c_{4123}(q) & = & -8(q-1)/(q+1)^4\\ 
        c_{4132}(q) & = & 16(q-1)^2/(q+1)^6\\ 
        c_{4213}(q) & = & 0\\       
        c_{4231}(q) & = & -32(q-1)/(q+1)^6\\       
        c_{4312}(q) & = & -32(q-1)/(q+1)^6\\       
        c_{4321}(q) & = & 64/(q+1)^6. \eeas
Although many values of $c_w(q)$ appear to be ``nice,'' not all are as
nice as the above data suggests. For instance,
  \beas c_{12453}(q) & = & 4(q^2+6q+1)(q-1)^4/(q+1)^8\\
    c_{13245}(q) & = & -2(q^4-8q^3-2q^2-8q+1)(q-1)^5/(q+1)^{10}\\
  c_{13425}(q) & = & -4(q^6-6q^5-33q^4+12q^3-33q^2-6q+1)(q-1)^2/
           (q+1)^{10}. \eeas 
We will prove two results about the $c_w(q)$'s.

\begin{theorem} \label{thm:qpromote}
Let $\mathrm{id}$ denote the identity permutation in $\sn$. Then
  $$ c_{\mathrm{id}}(q) =\left( \frac{q-1}{q+1}\right)^{\lfloor
    n/2\rfloor}. $$
\end{theorem}

\emph{Proof} (sketch). I am grateful to Monica Vazirani for
assistance with the following proof. Define a scalar product on $\hn$
by
   $$ \langle T_u,T_v\rangle = q^{\ell(u)}\delta_{uv}, $$
where $\ell(u)$ denotes the number of inversions of $u$ (i.e., the
length of $u$ as an element of the Coxeter group $\sn$). Then for any
$g,h\in\hn$ we have
  $$ \langle T_ig,h\rangle=\langle g,T_ih\rangle $$
and
  $$ \langle gT_i,h\rangle=\langle g,hT_i\rangle. $$
Since $E_i^2=1$ it follows that
  \beq \langle E_i gE_i,1 \rangle=\langle g,1\rangle.
      \label{eq:hsp} \eeq
Now 
  $$ c_{\mathrm{id}}(q) = \langle E_1E_2\cdots E_{n-1}E_1E_2\cdots
  E_{n-2}\cdots E_1E_2E_1,1\rangle. $$
Using equation~\eqref{eq:hsp} and the commutation relation $E_iE_j =
E_jE_i$ if $|i-j|\geq 2$, we obtain
  $$ c_{\mathrm{id}}(q) = \langle E_{n-1}E_{n-3}\cdots E_r,
     1\rangle, $$
where $r=1$ if $n$ is even, and $r=2$ if $n$ is odd. For any subset $S$
of $\{n-1,n-3,\dots,r\}$ we have
  $$ \prod_{i\in S}T_i = T_{\prod_{i\in S}s_i}. $$
(The $T_i$'s and $s_i$'s for $i\in S$ commute, so the above products
are well-defined.) Hence $\langle E_{n-1}E_{n-3}\cdots E_r,1\rangle$
is obtained by setting $T_i=0$ in each factor of the product $E_{n-1}
E_{n-2}\cdots E_r$, so we get
  $$ \langle E_{n-1}E_{n-3}\cdots E_r,1\rangle = \left(
    \frac{q-1}{q+1}\right)^{\lfloor n/2\rfloor}, $$
completing the proof.
\qed

\medskip
If $w=a_1 a_2\cdots a_n\in\sn$, then write $\widehat{w}$ for the reversal
$a_n\cdots a_2 a_1$. Equivalently, $\widehat{w}=w_0 w$, where
$w_0=n,n-1,\dots,1$ (the longest permutation in $\sn$). Our second
result on the polynomials $c_w(q)$ is the following.

\begin{theorem} \label{thm:qm1}
Let $w\in\sn$, and let $\kappa(w)$ denote the number of cycles of
$w$. Then $c_w(q)$ is divisible by $(q-1)^{n-\kappa(\widehat{w})}$.
\end{theorem}

\proof
Let $\bm{a}=(a_1,\dots,a_{\binom n2})$ be a reduced
decomposition of $w_0$. It is a well-known and simple consequence of
the strong exchange property for reduced decompositions (e.g.\
\cite[Thm.~1.4.3]{bj-br}) that if $k$ is the length of the longest
subsequence $(b_1,\dots,b_k)$ of $\bm{a}$ such that $s_{b_1} \cdots
s_{b_k}= w$, then $\binom n2-k$ is the minimum number of
transpositions $t_1,\dots, t_k$ for which $w=w_0 t_1\cdots t_k$. This
number is just $n-\kappa(w_0^{-1}w) = n-\kappa(w_0w)
=n-\kappa(\widehat{w})$, so $k=\binom n2-n+\kappa(\widehat{w})$. 

A fundamental property of the Hecke algebra $\hn$ is the
multiplication rule
 $$ T_uT_{k}=\left\{
   \begin{array}{ll}
  T_{u s_k}, & \mbox{if}\  l(u s_k)=l(u)+1, \\[.5em]
   q T_{u s_k}+(q-1)T_u, & \mbox{if}\  l(u s_k)=l(u)-1,
    \end{array} \right. $$
for any $u\in\sn$. Now consider the coefficient of $T_w$ in the
expansion of the product on the left-hand side of \eqref{eq:cwdef}.
For each factor $E_i=\frac{1}{q+1}(q-1-2T_i)$ we must choose a term
$(q-1)/(q+1)$ or $-2T_i/(q+1)$. If we choose $(q-1)/(q+1)$ then we
have introduced a factor of $q-1$. If we choose $-2T_i/(q+1)$ and
multiply some $T_u$ by it, then a $T_v$ so obtained satisfies either
$v=us_i$ or $v=u$; in the latter case a factor of $q-1$ is introduced.
It follows that every contribution to the coefficient of $T_w$ arises
from choosing a subsequence $(b_1,\dots,b_j)$ of the reduced
decomposition $(1,2,\dots,p-1,1,2,\dots,p-2,\dots,1,2,1)$ of $w_0$
such that $w=s_{b_1}\cdots s_{b_j}$, in which case we will obtain a
factor $(q-1)^{\binom n2-j}$. Since the largest possible value of $j$
is $\binom n2-n+\kappa(\widehat{w})$, it follows that $\binom n2-j\geq
n-\kappa(\widehat{w})$, completing the proof.  
\qed

\medskip
Theorem~\ref{thm:qm1} need not be best possible. For instance, some
values of $c_w(1)$ can be 0, such as $c_{2413}(q)$. For a nonzero
example, we have that $(q-1)^4$ divides $c_{2314}(q)$, but
$4-\kappa(4132)=2$. 
 

\end{document}